\documentclass[a4paper,12pt,table]{article}
\usepackage{color,float}
\usepackage{amsmath,amsfonts,amssymb,amscd}
\usepackage{enumerate}
\usepackage[ruled,linesnumbered]{algorithm2e}
\usepackage{indentfirst,graphicx, subfigure, epsfig,epstopdf,psfrag}

\input{epsf}
\usepackage{caption}

\newtheorem{thm}{Theorem}[section]

\newtheorem{lem}[thm]{Lemma}

\newtheorem{cor}[thm]{Corollary}

\def\qed{\hfill \nopagebreak\rule{5pt}{8pt}}

\title{\textbf{Hardness results for three kinds of colored connections of graphs\footnote{Supported by NSFC No.11871034 and 11531011.}}}
\author{\small  \ Zhong Huang, \ Xueliang Li \\[0.2cm]
\small Center for Combinatorics and LPMC \\
\small Nankai University, Tianjin 300071, China\\[0.2cm]
{\small  2120150001@mail.nankai.edu.cn, lxl@nankai.edu.cn }\\
}
\date{}
\begin{document}
\maketitle
\begin{abstract}
The concept of rainbow connection number of a graph was introduced by Chartrand et al. in 2008. Inspired by this concept, other concepts on
colored version of connectivity in graphs were introduced, such as the monochromatic connection number by Caro and Yuster in 2011, the proper connection number
by Borozan et al. in 2012, and the conflict-free connection number by Czap et al. in 2018, as well as some other variants of connection numbers later on.
Chakraborty et al. proved that to compute the rainbow connection number of a graph is NP-hard. For a long time, it has been tried to fix the computational
complexity for the monochromatic connection number, the proper connection number and the conflict-free connection number of a graph. However, it has not been
solved yet. Only the complexity results for the strong version, i.e., the strong proper connection number
and the strong conflict-free connection number, of these connection numbers were determined to be NP-hard.
In this paper, we prove that to compute each of the monochromatic connection number, the proper connection number and the conflict free connection number
for a graph is NP-hard. This solves a long standing problem in this field, asked in many talks of workshops and papers.


\end{abstract}

\section{Introduction}

All graphs considered in this paper are undirected, finite and simple. We follow \cite{BM} for notation and terminology not described here.
Let $G$ be a graph. We use $V(G), E(G), n(G), m(G), diam(G)$ to denote the vertex-set, edge-set, the number of vertices, number of edges,
and the diameter of $G$, respectively.

A path of an edge-colored graph $G$ is said to be a {\it rainbow path} if no two edges on the path have the same color.
The graph $G$ is called {\it rainbow connected} if every pair of distinct vertices of $G$ is connected by a rainbow path in $G$.
An edge-coloring of a connected graph is called a {\it rainbow connection coloring} if it makes the graph rainbow connected.
For a connected graph $G$, the {\it rainbow connection number $rc(G)$} of $G$ is defined as the smallest number of colors
that are needed in order to make $G$ rainbow connected. These concepts on the rainbow connection of graphs were introduced by
Chartrand et al. \cite{CJMZ} in 2008. The reader is referred to \cite{LSS, LS1, LS2} for details.

Inspired by the rainbow connection colorings and proper colorings in graphs, Borozan et al.~\cite{BFGMMMT} in 2012, and Andrews et al.~\cite{ALLZ}
independently, introduced the concept of proper-path coloring. Let $G$ be a nontrivial connected graph with an edge-coloring.
A path in $G$ is called a {\it proper path} if no two adjacent edges of the path receive the same color.
An edge-coloring $f$ of a connected graph $G$ is called a {\it proper-path coloring} if every pair of distinct vertices of $G$ are
connected by a proper path in $G$. If $k$ colors are used, then $f$ is called a {\it proper-path $k$-coloring}.
An edge-colored graph $G$ is {\it proper connected} if any two vertices of $G$ are connected by a proper path.
For a connected graph $G$, the minimum number of colors that are needed in order to make $G$ proper connected is
called the {\it proper connection number of $G$}, denoted by $pc(G)$. For more details, we refer to \cite{GLQ, HLQM, LLZ, DMP}
and a dynamic survey \cite{LM}.

Czap et al. \cite{CJV} in 2018 introduced the concept of conflict-free connection of graphs.
An edge-colored graph $G$ is called {\it conflict-free connected} if each pair of distinct
vertices is connected by a path which contains at least one color used on exactly one of its edges.
This path is called a {\it conflict-free path}, and this coloring is called a {\it conflict-free connection coloring} of $G$.
The {\it conflict-free connection number} of a connected graph $G$, denoted by $cfc(G)$,
is defined as the smallest number of colors needed to color the edges of $G$ such that $G$ is conflict-free connected.
More results can be found in \cite{DLLMZ, C1, C2}.

Note that a rainbow path has the maximum number of colors allowed on the path. Inspired by the concept of rainbow connection,
we would like to ask what about a path uses the minimum number of colors ? This will be the minimum cost when we try to establish
a point to point communication. As mentioned by Caro and Yuster \cite{CY}, it is a natural opposite counterpart of the rainbow connection colorings.
At first, we introduce the notion color-distance in the following.

Let $G$ be a nontrivial connected graph associated with an edge-coloring $f : E(G)\rightarrow \{1, 2, \ldots, p\}$, $p\in \mathbb{N}$.
Let $u$ and $v$ be two vertices of the edge-colored graph $G$. The {\it color-distance} between $u$ and $v$, denoted by $cd_f(u,v)$,
is the minimum number of distinct colors that a $uv$-path can have. It is easy to see that $cd_f(u,v)\geq 0$; $cd_f(u,v)=0$ if and only if $u=v$;
$cd_f(u,v)=cd_f(v,u)$; $cd_f(u,w)\leq cd_f(u,v)+cd_f(v,w)$. Then $(V, cd_f)$ is a metric space. We can also define color-distance between two sets
$U, V\subseteq V(G)$  of vertices: $cd_f(U,V)=max\{ cd_f(u,v): u\in U, v\in V\}$. Naturally, the {\it color-diameter} of $G$, denoted by $cD_f(G)$,
is the maximum color-distance among all pairs of vertices. We have $cD_f(G)=max\{cd_f(u,v): u,v\in V(G)\}$. The question is, when we fix the color-diameter,
how colorful can a coloring $f$ be ? That is to say, we already have $cD_f(G)\leq k$, and then we tend to maximize the
number of distinct colors $p$. We call $f$ a {\it $k$-color connection coloring} if $cD_f(G)\leq k$. The {\it $k$-color connection number} of a connected graph $G$,
denoted by $cc_k(G)$, is the maximum number of colors that are needed to make $f$ a $k$-color connection coloring. We call $f$ an
{\it extremal $k$-color connection coloring} if $f$ uses $cc_k(G)$ colors.

If $k=1$, then $cc_1(G)$ is exactly the {\it monochromatic connection number $mc(G)$} introduced by Caro and Yuster \cite{CY} in 2011.
More results on the monochromatic connection of graphs can be found in \cite{LW1, CLD, GGM, GLQZ, LW2}.

Now, for a given a graph $G$ we have got four parameters $rc(G)$, $pc(G)$, $cfc(G)$ and $mc(G)$. The question is that for a given graph, weather there exist
efficient ways to compute these parameters ? Chakraborty et al. \cite{SEAR} first solved this problem for the rainbow connection number $rc(G)$. They proved that computing $rc(G)$ is NP-hard.
However, for a long time it has been tried to fix the computational complexity for $pc(G)$, $cfc(G)$ and $mc(G)$ of a graph. Till now it has not been
solved, yet. Only the complexity results for the strong version, i.e., the strong proper connection number $spc(G)$ \cite{HY}
and the strong conflict-free connection number $scfc(G)$ \cite{JLZ}, were determined to be NP-hard.
In this paper, we prove that to compute each of $pc(G)$, $cfc(G)$ and $mc(G)$ of a graph is NP-hard.
This solves a long standing problem in this field, asked in many talks of workshops and papers; see \cite{LM, LMQ, LW1, C1, C2, JLZ}.

The paper is organized as follows. In Sections 2,3 and 4, we prove the NP-hardness for the proper connection number, conflict-free connection number
and monochromatic connection number, respectively. In the last section, Section 5, we work on the $k$-color connection number and give a linear-time algorithm
to compute the $k$-color connection number of trees.

\section{NP-hardnes for the proper connection number}

The classical 3-Satisfiability problem (3SAT for short) is well-known to be NP-complete.
Its variant Not-All-Equal 3 Satisfiability problem (NAE-3SAT for short) was also proved
to be NP-complete by Schaefer's dichotomy theorem \cite{S}. For ease of reading, we state
the variant problem as follows.

\noindent {\bf The NAE-3SAT problem:}

\noindent\textbf{INSTANCE:} A boolean formula $\phi$ in 3 conjunctive normal form (3CNF for short)
such that each clause is made up of three distinct literals.

\noindent\textbf{QUESTION:} Is there a satisfiable assignment for $\phi$ such that in each clause,
at least one literal is equal to true and at least one literal is equal to false ?

Our aim is to reduce the NAE-3SAT problem to the problem of deciding if $pc(G)=2$ for a graph $G$.
We then obtain our main result of this section.

\begin{thm}\label{pcthm}
For a graph $G$, deciding if $pc(G)=2$ is NP-complete. In particular, computing $pc(G)$ is NP-hard.
\end{thm}

\begin{figure}[htbp]
\vspace{-2cm}
\centering
\includegraphics[width=0.9\textwidth]{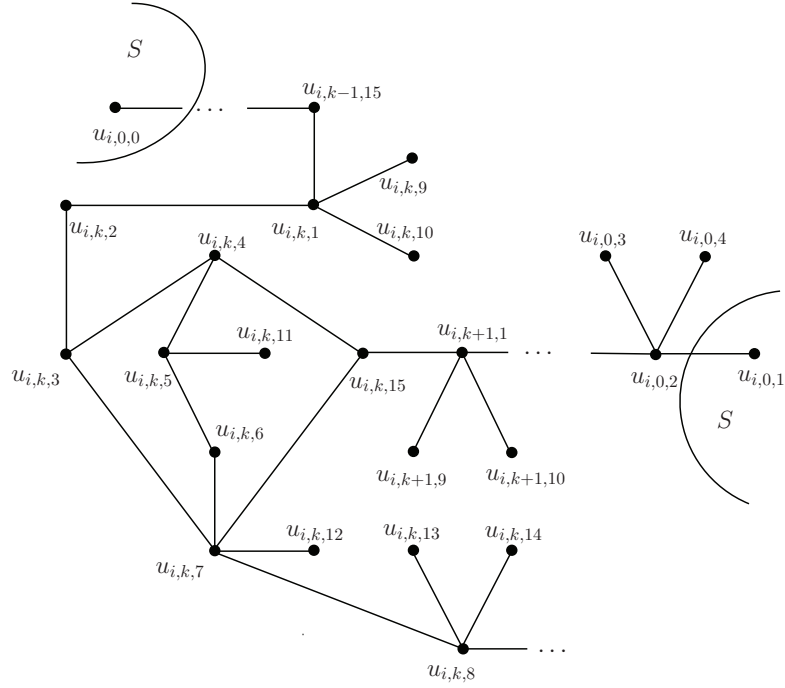}
\vspace{-6cm}
\caption{The variable garget $H_i$ corresponding to a variable $x_i$.}
\label{pc1}
\end{figure}

\begin{figure}[htbp]
\vspace{-2cm}
\centering
\includegraphics[width=0.9\textwidth]{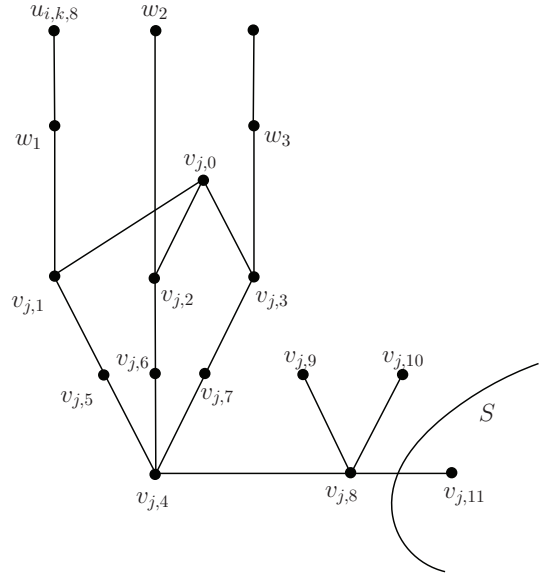}
\vspace{-8cm}
\caption{The graph $I_j$ representing for the clause $c_j$ such that the literals $l_{j,1}$, $l_{j,3}$ are positive and the literal $l_{j,2}$ is negative.}
\label{pc2}
\end{figure}

\begin{pf}
At first, we show that deciding if $pc(G)= 2$ is in NP. Actually,
Edmonds and Manoussakis showed that to check if a 2-edge-coloring of a graph makes
the graph proper connected is polynomial-time solvable; see \cite{LMQ}.
From their result, one can see that given a yes-instance of a 2-edge-colored connected
graph, it is polynomial-time checkable if the graph is proper connected.
Then we prove that the problem is NP-complete by reducing the NAE-3SAT problem to it.

Given a 3CNF $\phi=\wedge^m_{j=1}c_j$ on variables $x_1, x_2, \cdots, x_n$,
the clause $c_j$ is made up of three distinct literals $l_{j,1}, l_{j,2}, l_{j,3}$
for each $j\in [m]$. We construct a graph $G_\phi$ as follows:

First, we construct the variable garget $H_i$ for a variable $x_i$ for each $i\in [n]$
(see Figure \ref{pc1}). Assume that the variable $x_i$ appears in $m_i$ clauses of $\phi$.
For each $k\in[m_i]$, we start from a path $u_{i,k,1}u_{i,k,2}u_{i,k,3}u_{i,k,4}u_{i,k,5}u_{i,k,6}u_{i,k,7}u_{i,k,8}$.
Next we add pendent vertices as follows: each of the pendent vertices $u_{i,k,9}$ and $u_{i,k,10}$ is adjacent to $u_{i,k,1}$,
the pendent vertex $u_{i,k,11}$ is adjacent to $u_{i,k,5}$, the pendent vertex $u_{i,k,12}$ is adjacent to $u_{i,k,7}$,
and each of the pendent vertices $u_{i,k,13}$ and $u_{i,k,14}$ is adjacent to $u_{i,k,8}$.
Then we add a new vertex $u_{i,k,15}$, and a set of edges $\{u_{i,k,15}u_{i,k,4}, u_{i,k,15}u_{i,k,7}, u_{i,k,7}u_{i,k,3} \}$.
Now we get the graph $H_{i,k}$. The variable garget $H_i$ is obtained by adding a set of edges $\{u_{i,k,15}u_{i,k+1,1} | 1\leq k\leq m_i-1\}$
to the union of the graphs $H_{i,1}, H_{i,2}, \cdots ,H_{i,m_i}$.

Second, we construct the graph $I_j$ represent for the clause $c_j$ for each $j\in[m]$, see Figure \ref{pc2}.
We start from the vertices $v_{j,0}$ and $v_{j,4}$ such that there are three edge-disjoint paths
$v_{j,0}v_{j,1}v_{j,5}v_{j,4}$, $v_{j,0}v_{j,2}v_{j,6}v_{j,4}$ and $v_{j,0}v_{j,3}v_{j,7}v_{j,4}$ between them.
Next we add a new vertex $v_{j,8}$ such that $v_{j,8}$ is adjacent to $v_{j,4}$.
Then we add two pendent vertices $v_{j,9}$ and $v_{j,10}$ such that each of $v_{j,9}$ and $v_{j,10}$ is adjacent to $v_{j,8}$.
Then we add three vertices $w_1,w_2,w_3$ such that $w_1$ is adjacent to $v_{j,1}$, $w_2$ is adjacent to $v_{j,2}$,
and $w_3$ is adjacent to $v_{j,3}$. If a variable $x_i$ is positive such that $l_{j,1}=x_i$ in the clause $c_j$ for some $i\in[n]$,
then we add an edge $u_{i,k,8}w_1$ for some $k\in[m_i]$. If a variable $x_i$ is negative such that $l_{j,1}=\overline x_i$
in the clause $c_j$ for some $i\in[n]$, then we identity the vertices $u_{i,k,8}$ and $w_1$ for some $k\in[m_i]$.
We do this similarly for the literals $l_{j,2}$ and $l_{j,3}$.

Third, we add some extra vertices and edges to complete the construction of the graph $G_\phi$.
We construct a complete graph on the set of $2n+m+1$ vertices $S=\{u_{i,0,0}, u_{i,0,1},v_{j,11},s | i\in[n], j\in[m] \}$.
For each $i\in[n]$, we add a new vertex $u_{i,0,2}$ such that $u_{i,0,2}$ is adjacent to $u_{i, m_i, 15}$,
and add two pendent vertices $u_{i,0,3}$ and $v_{i,0,4}$ such that each of $u_{i,0,3}$ and $v_{i,0,4}$ is adjacent to $u_{i,0,2}$.
Next we add a set of edges $\{u_{i,0,0}u_{i,1,1}, u_{i,0,1}u_{i,0,2}, v_{j,8}v_{j,11}| i\in[n], j\in[m] \}$.
Now we get the graph $G_\phi$.

In the following, we only need to show that $pc(G_\phi)=2$ if and only if $\phi$ is satisfiable.

Suppose $pc(G_\phi)=2$. Let $f$ be a proper-path coloring of $G$ using 2 colors.

\textbf{Claim.} The three entries $u_{i,k-1,15}u_{i,k,1}$, $u_{i,k,15}u_{i,k+1,1}$ and $u_{i,k,7}u_{i,k,8}$ of $H_{i,j}$ share a same color for each $i\in[n], k\in[m_i]$.

\noindent\textbf{Proof of the Claim:}  Assume, to the contrary, that the three entries $u_{i,k-1,15}u_{i,k,1}$, $u_{i,k,15}u_{i,k+1,1}$ and $u_{i,k,7}u_{i,k,8}$ of $H_{i,j}$
do not share a same color. Notice that $u_{i,k,9}u_{i,k,1}u_{i,k,10}$ is the only path between the vertices $u_{i,k,9}$ and $u_{i,k,10}$ (see Figure \ref{pc1}).
Then we can assume that the edge $u_{i,k,9}u_{i,k,1}$ is colored with 1 and the edge $u_{i,k,10}u_{i,k,1}$ is colored with 2.
Similarly, we can assume that the edges $u_{i,k+1,9}u_{i,k+1,1}$, $u_{i,k+1,13}u_{i,k+1,8}$ are colored with 1 and the
edges $u_{i,k+1,10}u_{i,k+1,1}$, $u_{i,k+1,14}u_{i,k+1,8}$ are colored with 2. We distinguish the following two cases to complete the proof.

\textbf{Case 1.} The two entries $u_{i,k-1,15}u_{i,k,1}$ and $u_{i,k,15}u_{i,k+1,1}$ share a same color.

Without loss of generality, we can assume that the entries $u_{i,k-1,15}u_{i,k,1}$, $u_{i,k,15}u_{i,k+1,1}$ are colored with 1,
and the third entry $u_{i,k,7}u_{i,k,8}$ is colored with 2. Since there exists a proper path between vertices $u_{i,k,9}$ and $u_{i,k+1,10}$,
we can assume that the proper path is $u_{i,k,9}u_{i,k,1}u_{i,k,2}u_{i,k,3}u_{i,k,4}u_{i,k,5}u_{i,k,6}u_{i,k,7}u_{i,k,15}u_{i,k+1,1}u_{i,k+1,10}$.
Since there is a proper path between vertices $u_{i,k,12}$ and $u_{i,k,13}$, it follows that the edge $u_{i,k,12}u_{i,k,7}$ is colored with 1.
Notice that there exists a proper path between vertices $u_{i,k,11}$ and $u_{i,k,12}$, then the edge $u_{i,k,11}u_{i,k,5}$ is colored with 2.
There is a proper path between vertices $u_{i,k,11}$ and $u_{i,k,13}$, then the edge $u_{i,k,3}u_{i,k,7}$ is colored with 1.
Also there is a proper path between vertices $u_{i,k,9}$ and $u_{i,k,11}$, then the edge $u_{i,k,4}u_{i,k,15}$ is colored with 1.
However, there is no proper path between vertices $u_{i,k,13}$ and $u_{i,k,10}$. Thus, $f$ is not a proper-path coloring, a contradiction.

\textbf{Case 2.} The two entries $u_{i,k-1,15}u_{i,k,1}$ and $u_{i,k,15}u_{i,k+1,1}$ do not share a same color.

Without loss of generality, we assume that the entries $u_{i,k-1,15}u_{i,k,1}$ is colored with 1,
and the other two entries $u_{i,k,7}u_{i,k,8}$, $u_{i,k,15}u_{i,k+1,1}$ are colored with 2.
Then the proper path between vertices $u_{i,k,9}$ and $u_{i,k,13}$ is
$u_{i,k,9}u_{i,k,1}u_{i,k,2}u_{i,k,3}u_{i,k,4}u_{i,k,5}u_{i,k,6}u_{i,k,7}\\u_{i,k,8}u_{i,k,13}$.
Since there is a proper path between vertices $u_{i,k,12}$ and $u_{i,k,13}$, it follows that the
edge $u_{i,k,12}u_{i,k,7}$ is colored with 1. Notice that there exists a proper path between vertices
$u_{i,k,11}$ and $u_{i,k,12}$, then the edge $u_{i,k,11}u_{i,k,5}$ is colored with 2.
Suppose that the edge $u_{i,k,4}u_{i,k,15}$ is colored with 2. Then the edge $u_{i,k,7}u_{i,k,15}$ is colored with 1.
There is no proper path between vertices $u_{i,k,12}$ and $u_{i,k+1,9}$, then $f$ is not a proper-path coloring, a contradiction.
We have that the edge $u_{i,k,4}u_{i,k,15}$ is colored with 1. Since there is a proper path between vertices $u_{i,k,11}$ and $u_{i,k+1,9}$,
the edge $u_{i,k,7}u_{i,k,15}$ is colored with 1. Also, there is a proper path between vertices $u_{i,k,10}$ and $u_{i,k,11}$,
then the edge $u_{i,k,3}u_{i,k,7}$ is colored with 1. However, there is no proper path between vertices $u_{i,k,12}$ and $u_{i,k,14}$.
Thus, $f$ is not a proper-path coloring, a contradiction.

As a result, the three entries $u_{i,k-1,15}u_{i,k,1}$, $u_{i,k,15}u_{i,k+1,1}$ and $u_{i,k,7}u_{i,k,8}$ of $H_{i,j}$ share a same color for each $i\in[n], k\in[m_i]$. \qed

Now we have that the edges in the edge set $\{u_{i,k,7}u_{i,k,8} | k\in[m_i]\}$ share a same color for each $i\in[n]$.
If the edge $u_{i,0,0}u_{i,1,1}$ is colored with 1, then we set $x_i=0$ for each $i\in[n]$. If the edge $u_{i,0,0}u_{i,1,1}$ is colored with 2,
then we set $x_i=1$ for each $i\in[n]$. Let $s,t\in \{v_{j,1},v_{j,2},v_{j,3},v_{j,4}\}$ be two vertices,
$P_{s,t}$ be a path from $s$ to $t$ such that $P_{s,t}$ is in a graph $I_j$ for some $j\in[m]$. Then the length of $P_{s,t}$ is even.
Thus, two of the four edges $v_{j,1}w_1, v_{j,2}w_2, v_{j,3}w_3, v_{j,4}v_{j,8}$ should be colored with 1,
and the other two edges should be colored with 2 for each $j\in[m]$. Then $f$ admits a Not-All-Equal satisfiable assignment for the clause $c_j$.
As a result, we get a satisfiable assignment for $\phi$.

Suppose that $\phi$ is satisfiable. Given a Not-All-Equal satisfiable assignment on the variables $x_1,x_2, \cdots x_n$,
we color the graph $G$ as follows:

If $x_i=0$, then $u_{i,0,0}u_{i,1,1}$ is colored with 1. For each $k\in[m_i]$, the edges
$u_{i,k,9}u_{i,k,1}$, $u_{i,k,2}u_{i,k,3}$, $u_{i,k,4}u_{i,k,5}$, $u_{i,k,6}u_{i,k,7}$, $u_{i,k,8}u_{i,k,13}$, $u_{i,k,15}u_{i,k+1,1}$
are colored with 1, and the remaining edges in the graph $H_{i,k}$ is colored with 2. If $x_i=1$, then $u_{i,0,0}u_{i,1,1}$ is colored with 2.
For each $k\in[m_i]$, the edges $u_{i,k,9}u_{i,k,1}$, $u_{i,k,2}u_{i,k,3}$, $u_{i,k,4}u_{i,k,5}$, $u_{i,k,6}u_{i,k,7}$, $u_{i,k,8}u_{i,k,13}$,
$u_{i,k,15}u_{i,k+1,1}$ are colored with 2, and the remaining edges in the graph $H_{i,k}$ is colored with 1.
For each $i\in[n]$, the edge $u_{i,0,2}u_{i,0,3}$ shares a same color with the edge $u_{i,0,0}u_{i,1,1}$, and the
edges $u_{i,0,2}u_{i,0,4}$, $u_{i,0,2}u_{i,0,1}$ receive the distinct color. If a variable $x_i$ is
positive $x_i=l_{j,1}$ in the clause $c_j$ for some $k\in[m_i]$, then the edge $v_{j,1}w_1$ shares a same color with
the edge $u_{i,0,0}u_{i,1,1}$ and the edge $u_{i,k,8}w_1$ receives a distinct color. If a variable $x_i$ is negative
such that $\overline x_i=l_{j,1}$ in the clause $c_j$ for some $k\in[m_i]$, then the edge $v_{j,1}w_1$ and the edge $u_{i,0,0}u_{i,1,1}$
receive distinct colors. The coloring is similar for the literals $l_{j,2}$ and $l_{j,3}$. Without loss of generality,
we can assume that the edges $v_{j,1}w_1$ and $v_{j,2}w_2$ are colored with 1, the edge $v_{j,3}w_3$ is colored with 2, since $\phi$ is satisfiable.
Then the edges $v_{j,1}v_{j,0}$, $v_{j,1}v_{j,5}$, $v_{j,2}v_{j,0}$, $v_{j,2}v_{j,6}$, $v_{j,4}v_{j,7}$,
$v_{j,8}v_{j,10}$, $v_{j,8}v_{j,1}$ are colored with 1, the edges $v_{j,3}v_{j,0}$, $v_{j,3}v_{j,7}$,
$v_{j,4}v_{j,5}$, $v_{j,4}v_{j,6}$, $v_{j,4}v_{j,8}$, $v_{j,8}v_{j,9}$ are colored with 2. Let $p\neq s$ be a vertex in $S$.
Then there exists a unique vertex $q\in V(G_\phi)\setminus S$ such that $p$ is adjacent to $q$.
If the edge $pq$ is colored with 1, then we color the edge $ps$ with 2. If the edge $pq$ is colored with 2,
then we color the edge $ps$ with 1. Let $p_1,p_2\in S$ be two vertices. If the edges $p_1s$ and $p_2s$ share a same color 1,
then we color the edge $p_1p_2$ with 1. If the edges $p_1s$ and $p_1s$ share a same color 2, then we color the edge $p_1p_2$ with 2.
The remaining edges can be arbitrarily colored.

One can verify that there exists a proper path between every pair of vertices in the graph $G_\phi$. Thus, $pc(G_\phi)=2$.
The proof is now complete. \qed

\end{pf}

\section{NP-hardness for the conflict-free connection number}

Let $C(G)$ be the subgraph of $G$ induced by the set of cut-edges of $G$,
and let $h(G) = max\{cfc(T): \text{ $T$ is a component of $C(G)$}\}$. The following result
was obtained in \cite{CJV}.

\begin{lem}\cite{CJV}\label{cfcnp}
If $G$ is a connected graph with cut-edges, then
 $$h(G) \leq cfc(G) \leq h(G) + 1.$$ \\
Moreover, the bounds are sharp.
\end{lem}

In the following, our aim is to reduce the NAE-3SAT problem to the problem of deciding whether $cfc(G)=h(G)$ or $h(G)+1$ for the case when $h(G)=2$.
We then know that deciding whether $cfc(G)=2$ or $3$ is NP-complete by Theorem \ref{cfcthm}. At first, we need the following lemma.
\begin{lem}\label{cfclem}
Let $G$ be a graph such that $G$ contains two edge-disjoint Hamilton paths $P_1$ and $P_2$.
The edges of $P_1$ is colored with 1 and the edges of $P_2$ is colored with 2.
Then between any pair of vertices of $G$, there exists a monochromatic color 1 path,
a monochromatic color 2 path, a conflict-free path with color 1 used only once and a conflict-free path with color 2 used only once.
\end{lem}
\begin{pf}
Let $u$ and $v$ be two vertices of $G$. Clearly, $u$ and $v$ lie on $P_1$, then the path $uP_1v$ is a monochromatic color 1 path.
Similarly, the path $uP_2v$ is a monochromatic color 2 path.
Now consider the path $P_2$. Let $e$ be the edge incident with $u$ on the path $uP_2v$.
Then $P_2-e$ contains two component $P_u$ and $P_v$ ($u$ lies on $P_u$). Since $P_1$ is a Hamilton path,
there exists an edge $xy$ on $P_1$ such that $x\in V(P_u), y\in V(P_v)$. Now we have that the path $uP_2xyP_2v$ is a conflict-free path with color 1 used only once.
Similarly, there exists a conflict-free path with color 2 used only once.
\end{pf}

Next we need the lemma given by M Ji et al.~\cite{JLZ} as following:

\begin{lem} \label{jlz} \cite{JLZ}
Given a connected graph $G$ and a coloring $f$ : $E \ (\text{or} \ V )\rightarrow \{1, 2,\cdots,$ $k\}$ $(k \geq 1)$ of
$G$, determining whether or not $G$ is, respectively, conflict-free connected, conflict-free vertex-connected,
strongly conflict-free connected under the coloring $f$ can be done in polynomial-time.
\end{lem}

Our main result of this section is as follows.

\begin{figure}[htbp]
\vspace{-2cm}
\centering
\includegraphics[width=0.9\textwidth]{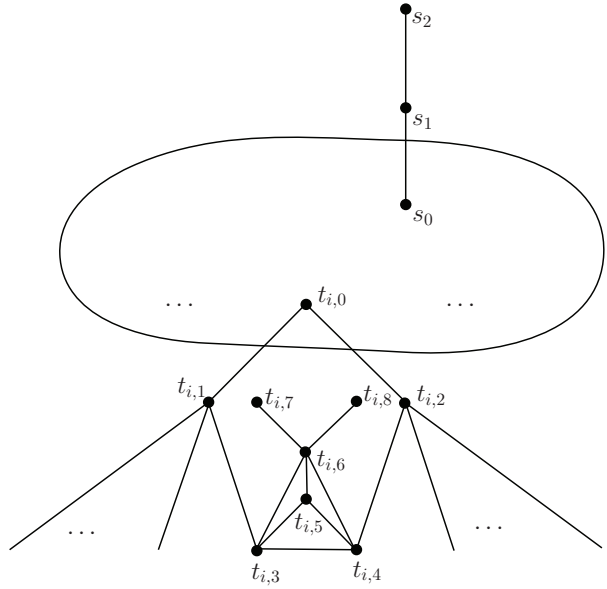}
\vspace{-8cm}
\caption{The variable garget $H_i$ corresponding to a variable $x_i$.}
\label{cfc1}
\end{figure}

\begin{figure}[htbp]
\vspace{-2cm}
\centering
\includegraphics[width=1.1\textwidth]{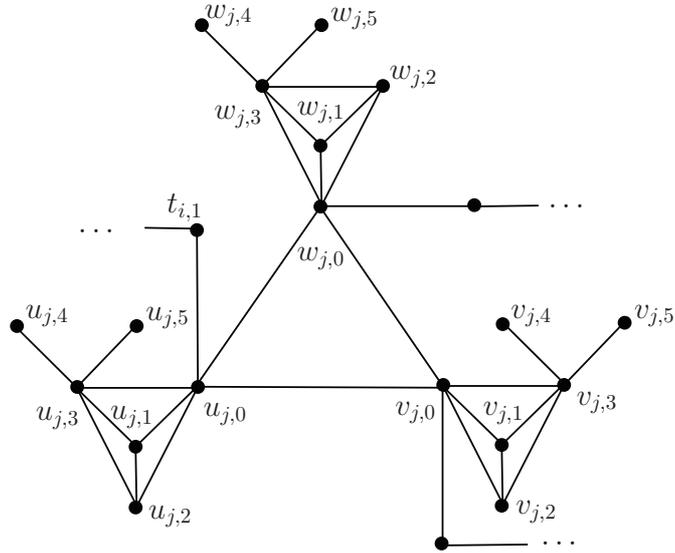}
\vspace{-12cm}
\caption{The graph $I_j$ representing for the clause $c_j$ such that the variable $x_i=l_{j,1}$ is positive in $c_j$.}
\label{cfc2}
\end{figure}

\begin{thm}\label{cfcthm}
Let $G$ be a graph such that $h(G)=2$. Then deciding whether $cfc(G)=2$ or 3 is NP-complete.
In particular, computing $cfc(G)$ is NP-hard.
\end{thm}

\begin{pf}
As a result of Lemma \ref{jlz}, the problem is in NP. Then we prove that the problem is NP-complete by reducing the NAE-3SAT problem to it.
Given a 3CNF $\phi=\wedge^m_{j=1}c_j$ on variables $x_1, x_2, \cdots, x_n$, the clause $c_j$ is made up of three distinct literals
$l_{j,1}, l_{j,2}, l_{j,3}$ for each $j\in [m]$. We construct a graph $G_\phi$ as follows.

First, we construct the variable garget $H_i$ for a variable $x_i$ for each $i\in [n]$.
We start from a complete graph on four vertices $\{t_{i,3}, t_{i,4}, t_{i,5}, t_{i,6}\}$.
Next we add two pendent vertices $t_{i,7}, t_{i,8}$ such that each of the vertices $t_{i,7}, t_{i,8}$ is adjacent to $t_{i,6}$.
Then we get $H_i$ by adding a cycle $t_{i,0}t_{i,1}t_{i,3}t_{i,4}t_{i,2}t_{i,0}$ of length 5.

Second, we construct the graph $I_j$ representing for the clause $c_j$ for each $j\in[m]$.
We start from a complete graph on four vertices $\{u_{j,0}, u_{j,1}, u_{j,2}, u_{j,3}\}$.
Next we add two pendent vertices $u_{j,4}$ and $u_{j,5}$ such that each of $u_{j,4}$ and $u_{j,5}$ is adjacent to $u_{j,3}$.
Let the resulting graph be $I_{j,u}$. Then we add two copies $I_{j,v}$ and $I_{j,w}$ of $I_{j,u}$. Finally,
we get $I_j$ by adding the set of edges $\{u_{j,0}v_{j,0}, v_{j,0}w_{j,0}, w_{j,0}x_{j,0}\}$ to the union of $I_{j,u}$, $I_{j,v}$ and $I_{j,w}$.
If a literal $l_{j,1}=x_i$ is positive in the clause $c_j$ for some $i\in[n]$, then $u_{j,0}$ is adjacent to $t_{i,1}$.
If a literal $l_{j,1}=\overline x_i$ is negative in the clause $c_j$ for some $i\in[n]$,
then $v_{j,0}$ is adjacent to $t_{j,2}$. We do this similarly for the literals $l_{j,2}$ and $l_{j,3}$.

Third, we add some extra vertices and edges to complete the construction of the graph $G_\phi$.
We construct a complete graph on $n+1$ vertices $\{t_{1,0}, t_{2,0}, \cdots, t_{n,0}, s_0\}$.
Next we add two vertices $s_1$ and $s_2$ such that $s_0s_1s_2$ is a path of length two. Then the resulting graph is $G_\phi$.

In the following, we only need to show that $cfc(G_\phi)=2$ if and only if $\phi$ is satisfiable.

Suppose $cfc(G_\phi)=2$. Notice that $u_{j,4}u_{j,3}u_{j,5}$ is the only path between the vertices $u_{j,4}$ and $u_{j,5}$ (see Figure \ref{cfc2}).
Then we can assume that the edge $u_{j,4}u_{j,3}$ is colored with 1 and the edge $u_{j,5}u_{j,3}$ is colored with 2.
Similarly, we can assume that the edges $v_{j,4}v_{j,3}$, $w_{j,4}w_{j,3}$ are colored with 1
and the edges $v_{j,5}v_{j,3}$, $w_{j,5}w_{j,3}$ are colored with 2. We can assume that $s_0s_1$ is colored with 1 and $s_1s_2$ is colored with 2,
since $s_0s_1s_2$ is the only path between vertices $s_0$ and $s_2$ (see Figure \ref{cfc1}).
There exists a monochromatic color 1 path between vertices $s_0$ and $u_{j,3}$, since color 1 has been used at least twice
on the conflict-free path between vertices $s_2$ and $u_{j,4}$. Similarly, there exists a monochromatic color 2 path
between $s_0$ and $u_{j,3}$. The graph $I_j$ has exact three entries, then these three edges can not have a same color;
otherwise, we can not get both a monochromatic color 1 path and a monochromatic color 2 path between $s_0$ and $u_{j,3}$.

\textbf{Claim.} There exists no monochromatic path going through some $I_j$.

\noindent\textbf{Proof of the Claim:} Suppose that there exists a monochromatic color 1 path going through some
two entries of $I_j$. Then the third entry should be colored with 2. There exists a monochromatic color 1 path
and a monochromatic color 2 path between vertices $s_0$ and $u_{j,0}$, since any path between $s_0$ and $u_{j,3}$ goes
through the vertex $u_{j,0}$. Similarly, there exists a monochromatic color 1 path and a monochromatic color 2 path between
vertices $s_0$ and $v_{j,0}$; there exists a monochromatic color 1 path and a monochromatic color 2 path between vertices $s_0$ and $w_{j,0}$.
Then we have that the edge $u_{j,0}v_{j,0}$ is colored with 2. Next, we color the edge $u_{j,0}w_{j,0}$ with 1 and the edge $v_{j,0}w_{j,0}$ with 2.
Hence, there exists no monochromatic color 1 path going through $I_j$ for any $j\in[m]$. Similarly, there exists no monochromatic color 2 path going through some $I_j$. \qed

Notice that $t_{i,7}t_{i,6}t_{i,8}$ is the only path between vertices $t_{i,7}$ and $t_{i,8}$ (see Figure \ref{cfc2}).
Then we can assume that the edge $t_{i,7}t_{i,6}$ is colored with 1 and the edge $t_{i,8}t_{i,6}$ is colored with 2.
We have that there exists a monochromatic color 1 path and a monochromatic color 2 path between vertices $s_0$ and $t_{i,6}$.
Then the two edges $t_{i,0}t_{i,1}$ and $t_{i,0}t_{i,2}$ are colored with distinct colors by our Claim.
If $t_{i,1}$ ($t_{i,2}$, respectively) is incident with the entry $u_{j,0}t_{i,1}$ ($u_{j,0}t_{i,2}$, respectively) of $I_j$,
then the edge $u_{j,0}t_{i,1}$ ($u_{j,0}t_{i,2}$, respectively) and the edge $t_{i,0}t_{i,1}$ ($t_{i,0}t_{i,2}$, respectively)
share a same color to ensure the monochromatic paths between vertices $s_0$ and $u_{j,0}$. Similar thing happens for the remaining two entries of $I_j$.
If $t_{i,0}t_{i,1}$ is colored with 1, then we set $x_i=0$. If $t_{i,0}t_{i,1}$ is colored with 2, then we set $x_i=1$. As a result,
the variables admit a Not-All-Equal satisfiable assignment of $\phi$.

Suppose that $\phi$ is satisfiable. Given a Not-All-Equal satisfiable assignment on variables $x_1,x_2, \cdots x_n$, \
we color the graph $G$ as follows.

If $x_i=0$, then the edges incident with the vertex $t_{i,1}$ are colored with 1 and the edges incident with the
vertex $t_{i,2}$ is colored with 2 for each $i\in[n]$. If $x_i=1$, then the edges incident with the vertex $t_{i,1}$
are colored with 2 and the edges incident with the vertex $t_{i,2}$ is colored with 1 for each $i\in[n]$.
The three entries of $I_j$ do not share a same color, since $\phi$ is satisfiable for each $j\in[m]$.
Without loss of generality, assume that the entry incident with $u_{j,0}$ and the entry incident with $v_{j,0}$ are colored with 1,
and the entry incident with $w_{j,0}$ is colored with 2 for each $j\in[m]$. Then we color the two edges $u_{j,0}v_{j,0}$ and $u_{j,0}w_{j,0}$ with 2,
and color the edge $v_{j,0}w_{j,0}$ with 1 for each $j\in[m]$. Next the edges $u_{j,4}u_{j,3}$, $v_{j,4}v_{j,3}$, $w_{j,4}w_{j,3}$ are colored with 1,
and the edges $u_{j,5}u_{j,3}$, $v_{j,5}v_{j,3}$, $w_{j,5}w_{j,3}$ are colored with 2 for each $j\in[m]$.
The edge $t_{i_7}t_{i,6}$ is colored with 1 and the edge $t_{i,8}t_{i,6}$ is colored with 2 for each $i\in[n]$.
The remaining parts consist of complete subgraphs of order at least 4. By induction, we know that a complete graph of order at least 4
contains two edge-disjoint Hamilton paths. For each complete subgraph, one Hamilton path is colored with 1, and the other is colored with 2.

Next we should verify that there exists a conflict-free path between every pair of vertices.
Notice that there exist both a monochromatic color 1 path and a monochromatic color 2 path
between $s_0$ and $u_{j,0}$ for each $j\in[m]$. Let $y, z \in \cup_{i=1}^n V(H_i)\cup \{s_0, s_1, s_2\}$ be two vertices.
Then there exists a conflict-free path between $y$ and $z$ by Lemma \ref{cfclem}.
Let $y\in \cup_{i=1}^n V(H_i)\cup \{s_0, s_1, s_2\}, z\in V(I_j)$ be two vertices for some $j\in[m]$.
Then there exists a conflict-free path between $y$ and $z$ by Lemma \ref{cfclem}.
Let $y, z\in V(I_j)$ be two vertices for some $j\in [m]$. Then there exists a conflict-free path between $y$ and $z$ by Lemma \ref{cfclem}.
Let $y\in V(I_{j_1})\setminus \{u_{j_1,0}, v_{j_1,0}, w_{j_1,0}\}$, $z\in V(I_{j_2})$ be two vertices for some $j_1,j_2\in [m]$.
Then there exists a conflict-free path between $y$ and $z$ by Lemma \ref{cfclem}.
Let $y\in \{u_{j_1,0}, v_{j_1,0}, w_{j_1,0}\}$, $z\in V(I_{j_2})\setminus \{u_{j_2,0}, v_{j_2,0}, w_{j_2,0}\}$ be two vertices for some $j_1,j_2\in [m]$.
Then there exists a conflict-free path between $y$ and $z$ by Lemma \ref{cfclem}.
Without loss of generality, for some $j_1,j_2\in[m]$, assume that the entry incident with $u_{j_1,0}$ and the entry incident with $v_{j_1,0}$ are
colored with 1, and the entry incident with $w_{j_1,0}$ is colored with 2.
Then we color the two edges $u_{j_1,0}v_{j_1,0}$ and $u_{j_1,0}w_{j_1,0}$ with 2, and color the edge $v_{j_1,0}w_{j_1,0}$ with 1.
Let $y=u_{j_1,0}$, $z\in \{u_{j_2,0}, v_{j_2,0}, w_{j_2,0}\}$ be two vertices.
Then we can find a conflict-free path $P_{y,z}$ from $y$ to $z$ with color 2 which is the unique color such that edge $u_{j_1,0}v_{j_1,0}$ is in $E(P_{y,z})$.
Let $y=w_{j_1,0}$, $z\in \{u_{j_2,0}, v_{j_2,0}, w_{j_2,0}\}$ be two vertices.
Then we can find a conflict-free path $P_{y,z}$ from $y$ to $z$ with color 2 which is the unique color such that edge $u_{j_1,0}w_{j_1,0}$ is in $E(P_{y,z})$.
Let $y=v_{j_1,0}$, $z\in \{u_{j_2,0}, v_{j_2,0}, w_{j_2,0}\}$ be two vertices.
Then we can find a conflict-free path $P_{y,z}$ from $y$ to $z$ with color 1 which is the unique color such that edge $v_{j_1,0}w_{j_1,0}$ is in $E(P_{y,z})$.
As a result, we have that $cfc(G)=2$. The proof is now complete. \qed

\end{pf}

\section{NP-hardness for the monochromatic connection number}

In the last section we will deal with the complexity for the monochromatic connection number. The following result
was obtained in \cite{CY}.

\begin{thm}{\upshape \cite{CY}}\label{mc=m-n+2}
Let $G$ be a connected graph with $n > 3$ vertices and $m$ edges. If $G$
satisfies one of the following conditions, then $mc(G) = m - n + 2$.
\begin{enumerate}
\item $G$ is 4-connected.
\item $G$ is triangle-free.
\item $\Delta(G) < n- \frac{2m-3(n-1)}{n-3}$. In particular, this holds if $\Delta(G) \leq \frac{n+1}{2}$ or $\Delta(G) \leq n - \frac{2m}{2n}$.
\item $diam(G) \geq 3$.
\item $G$ has a cut vertex.
\end{enumerate}
\end{thm}

Inspired by the concept of a monochromatic path which uses exact one color,
we consider a path using at most $k$ distinct colors. Next, we will use the color-distance as a tool.
Given a positive integer $k$ and a graph $G$, if one asks for an edge-coloring of $G$ satisfying
that for any pair of vertices of $G$ there exists a path in $G$ which uses at most $k$ distinct colors between the two vertices,
then the maximum number of distinct colors for the edge-coloring $G$ is $cc_k(G)$.
We can consider $cc_k(G)$ in an opposite way which counts the total number of ''wasted" colors in the following lemma.

\begin{lem}\label{Waste}
Let $f$ be an extremal $k$-color connection coloring of $G$. Then each color forms a forest.
For a color $r$, denote by $F_r$ the subgraph induced by the set of edges colored with $r$.
We call $F_r$ the {\it color forest} of the color $r$. A color forest with $m(F_r)$
edges is said to {\it waste} $m(F_r)-1$ colors. Hence, $F_r$ is said to be nontrivial if and only if $m(F_r)-1>0$.
Let $W_k(G)$ be the total number of wasted colors. Then $cc_k(G)=m(G)-W_k(G)$.
\end{lem}
\begin{pf}
We claim that the subgraph $F_r$ induced by the set of edges colored with $r$ is a forest. Assume to the contrary,
there is a cycle in $F_r$. Then it is possible to choose any edge in this cycle and color it with a fresh color
while still maintaining a $k$-color connecting coloring, a contradiction.

Next, the $k$-color connection number counts the number of distinct color forests.
Hence, we have $cc_k(G)=\sum_c 1 =m(G)- \sum_c (m(F_c)-1)= m(G)-W_k(G)$. \qed
\end{pf}

\begin{lem}\label{range}
For a connected graph $G$, $0\leq W_k(G)\leq n-1-k$.
\end{lem}
\begin{pf}
By definition, we know that $W_k(G)=0$ if and only if $diam(G)$ is at most $k$.
Let $H$ be a connected induced subgraph of $G$. We can construct a $k$-color connecting coloring
of $G$ by adding $m(G)-m(h)$ new colors to a $k$-color connecting coloring of $H$,
and hence $cc_k(G)\geq cc_k(H)+m(G)-m(H)$. Then we have $W_k(G)\leq W_k(H)$ by Lemma \ref{Waste}.

Next we consider the upper bound of $W_k(G)$. Let $T$ be a spanning tree of graph $G$.
Then $cc_k(T)$ can be calculated in Section 4. Hence, $m(G)-cc_k(T)$ is an alternative upper bound.
In general, we have $W_k(G)\leq n-1-k$ by our Algorithms \ref{odd} and \ref{even} in next section.
Moreover, $W_k(G)$ achieves its upper bound if and only if $G$ is a path. \qed
\end{pf}

\begin{lem}\label{NP}
Let $k$ be a positive integer. Given a connected graph $G$ associated with an edge-coloring $f$ on
$n$ vertices, one can decide if $cD(G)\leq k$ in polynomial time.
\end{lem}
\begin{pf}
Let $[l]$ be the color set of $f$. Clearly, $l\leq m(G)$. Next we choose a pair of vertices $v,t$ in $G$ arbitrarily.
Let $C$ be a color set chosen from $[l]$ such that $|C|=k$. In this step, we have exact $\tbinom{l}{k}$ choices.
Since $k$ is a constant and $m(G)\leq n^2$, it follows that we can check each choice in $\mathcal{O}(n^{2k})$ time.
Then we star form a vertex set $S_0=\{v\}$.
If there exist two vertices $u\in S_i$ and $u^\prime\not\in S_i$ such that the edge $uu^\prime$ receives a color $c\in C$,
then $S_{i+1}=S_i\cup u^\prime$ for $i\in[n]$. If there do not exist two vertices $u\in S_i$ and $u^\prime\not\in S_i$ such that the edge $uu^\prime$ receives a color $c\in C$,
then $S_{i+1}=S_i$ for $i\leq n(G)$. In this step, we should check the neighbour of vertex set $S_i$ for each $i\in[n]$. And this can be done in $\mathcal{O}(n^2)$ time.
Since there are at most $\tbinom{n}{2}$ pairs of vertices, one can decide if $cD(G)\leq k$ in $\mathcal{O}(n^{2k+4})$ time. \qed
\end{pf}

\begin{figure}[htbp]
\centering
\includegraphics[width=0.5\textwidth]{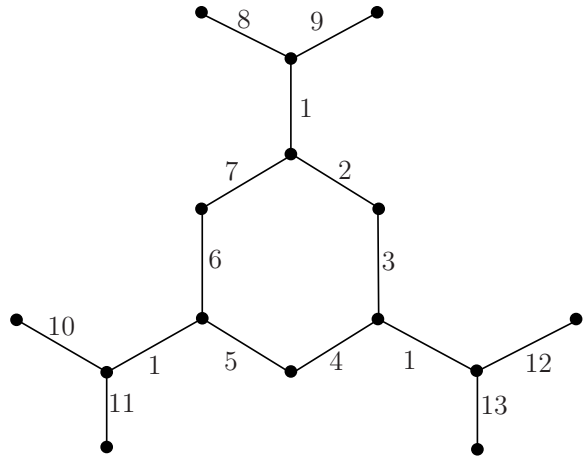}
\caption{$cc_5(G)=13$ and the color forest $F_1$ consists of 3 components.}
\label{mc1}
\end{figure}

Let $f$ be an extremal $1$-color connection coloring of $G$, then the color forest $F_c$ is a tree,
and the tree is called a color tree $T_c$ in \cite{CY}. However $F_c$ is not always a tree for $k>1$, see Figure \ref{mc1}.

Next we consider the case when $k=1$.

Let $u,v$ be two vertices of $G$. If $u,v$ are adjacent, then $cd(u,v)=1$. If $u,v$ are nonadjacent, then we should waste colors to
construct a path that uses only one color between $u$ and $v$. Hence we consider the complement graph $\overline G$.
We say that a color tree $T_c$ connects the edge $uv\in \overline G$ when $u$ and $v$ lie in a same color tree.

Let $w$ be a vertex of degree $p$ in $\overline G$.
Then there are more than $p$ wasted edges related to $w$. We display the fact in the following lemma.

\begin{lem}\label{p}
Let $f$ be an extremal 1-color connection coloring of $G$, and the degree of a vertex $w\in \overline G$
be $p$. Then $\sum_{T_c\in S} (m(T_c)-1)\geq p$, in which $S=\{T_c: w\in V(T_c), c\in [l]\}$.
\end{lem}

\begin{pf}
Let $N[w]$ be the closed neighborhood of vertex $w \in \overline G$,
and $T_c$ be a nontrivial color tree such that it connects some edges
which are adjacent to $w$ in $\overline G$. Clearly, $w\in V(T_c)$ and we have $V(T_c)\nsubseteq N[w]$.
Suppose $V(T_c)\cap N[w]=t$. Then $T_c$ connects $t-1$ edges which are adjacent to $w$.
Since $V(T_c)\geq t+1$, $m(T_c)-1\geq t-1$. Summing up the inequality, and we have $\sum_{T_c\in S} (m(T_c)-1)\geq p$. \qed
\end{pf}

\begin{figure}[htbp]
\begin{center}
\mbox{\subfigure[$P_4$]{\includegraphics[height=30mm]{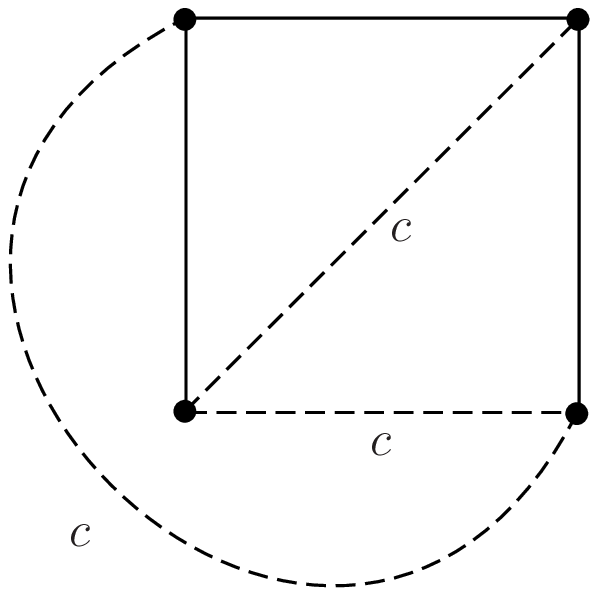}}\quad
\subfigure[$P_6$]{\includegraphics[height=50mm]{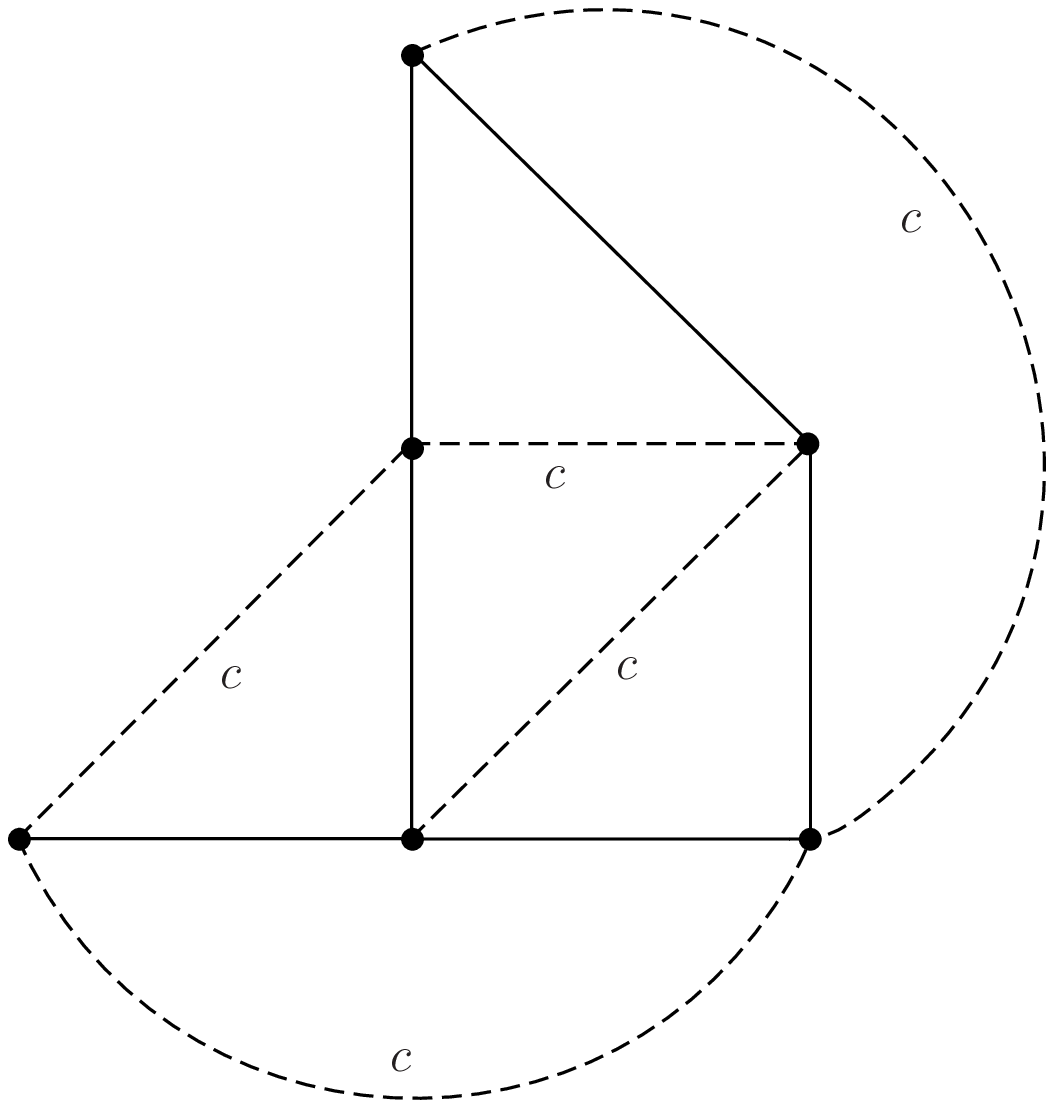}}
}
\caption{The nontrivial color tree $T_c$ is a path on four vertices or six vertices.}
\label{mc23}
\end{center}
\end{figure}

\begin{figure}[htbp]
\vspace{-2cm}
\centering
\includegraphics[width=1\textwidth]{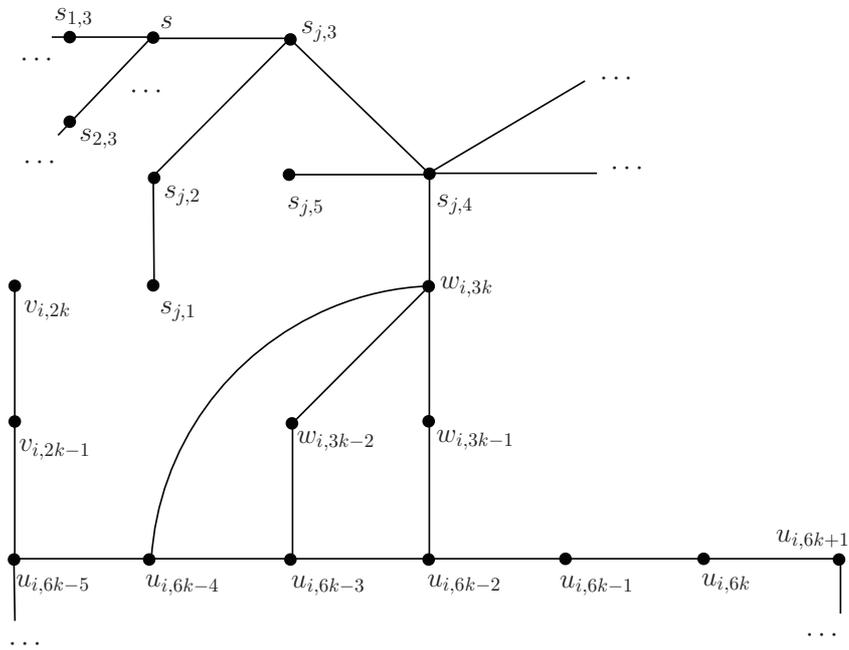}
\vspace{-10cm}
\caption{The variable $x_i$ is positive in clause $c_j$ for some $j\in[m]$.}
\label{mc4}
\end{figure}

Given a connected graph $G$, it is NP-complete to decide if $cc_1(G)$ attains its trivial lower bound $m(G)-n(G)+2$.

\begin{thm}\label{cck1}
Given a graph $G$, deciding whetehr $cc_1(G)= m(G)-n(G)+2$ or $cc_1(G)> m(G)-n(G)+2$ is NP-complete.
\end{thm}
\begin{pf}
The problem is in NP by Lemma \ref{NP}. We will prove that it is NP-complete by reducing the 3SAT problem to it.
Given a 3CNF $\phi=\wedge^m_{j=1}c_j$ on variables $x_1,x_2,\cdots,x_n$. The variable $x_i$ appears in $m_i$ clauses for each $i\in[n]$.
Notice that $m=\frac{1}{3}\Sigma_{i=1}^n m_i$. Then we construct a graph $G_{\phi}$ on $38m+1$ vertices as follows.

First, we construct the variable garget $H_i$ for a variable $x_i$, where $H_i$ is obtained from a cycle of order $6m_i$
with vertices $u_{i,1},u_{i,2},\cdots,u_{i,6m_i}$ such that for each $j\in[m_i]$, a new path $u_{i,6j-5}v_{i,2j-1}v_{i,2j}$ of length 2 is added to each vertex of $u_{i,6j}$.

Second, we construct the graph $I_j$ representing for the clause $c_j$ for $j\in[m]$. We start form a path $s_{j,1}s_{j,2}s_{j,3}s_{j,4}s_{j,5}$ of length 4.
Then we add an edge $s_{j,3}s$ such that $s_{j,3}$ is adjacent to the global vertex $s$ for each $k\in[m]$. If a variable $x_i$ is positive in clause $c_j$,
then we add three new vertices $w_{i,3k-2},w_{i,3k-1},w_{i,3k}$ and a set of edges
$\{w_{i,3k}s_{j,4}, w_{i,3k}w_{i,3k-1}, w_{i,3k}w_{i,3k-2}, w_{i,3k}u_{i,6k-1}, w_{i,3k-1}u_{i,6k-2}, w_{i,3k-2}u_{i,6k-3}\}$ for each $k\in[m_i]$.
If a variable $x_i$ is negative in clause $c_j$, then we add three new vertices $w_{i,3k-2},w_{i,3k-1},w_{i,3k}$ and a set of edges
$\{w_{i,3k}s_{j,4}, w_{i,3k}w_{i,3k-1}, w_{i,3k}w_{i,3k-2},\\ w_{i,3k}u_{i,6k-5}, w_{i,3k-1}u_{i,6k-4}, w_{i,3k-2}u_{i,6k-3}\}$ for each $k\in[m_i]$.

Third, let $V_0=\{v_{i,2k}, u_{i,6k-2}, s_{j,1}| j\in[m_i], i\in[n], k\in[m_i]\}$, $|V_0|=7m$.
Then we add the edges from an edge set $\{st | s,t\in V_0\}$. After deleting the edges of a Hamilton path $P$ on the vertex set $V_0$,
we get the graph $\overline {G_\phi}$. Now $G_\phi$ is defined.

In the following, we only need to show that $cc_1(G_{\phi})>m(G_\phi)-38m+1$ if and only if $\phi$ is satisfiable.

Suppose that $\phi$ is satisfiable. Then we define a 1-color-connection coloring $f$ on $G\phi$ using $m(G_\phi)-38l+2$ colors.
First, the edges on the path $P$ receive a same color, then the color tree $P$ wastes $7l-2$ colors.
Next, suppose $x_i=0$. Then we give a clockwise coloring such that path $u_{i,6k-5}v_{i,2k}u_{i,6k-4}v_{i,2k-1}$ is a color tree.
If a variable $x_i$ is positive in clause $c_j$ for some $j\in[m]$, then the path $u_{i,6k-3}w_{i,3k-1}u_{i,6k-4}u_{i,6k-2}$ is a color tree,
the path $w_{i,3k}u_{i,6k-1}w_{i,3k-1}w_{i,3k-2}$ is a color tree, and the path $u_{i,6k-1}u_{i,6k+1}u_{i,6k-2}u_{i,6k}$ is a color tree for some $k\in[m_i]$.
If a variable $x_i$ is negative in clause $c_j$ for some $j\in[m]$, then the path $w_{i,3k}u_{i,6k-2}u_{i,6k-4}w_{i,3k-2}$ $w_{i,3k-1}u_{i,6k-3}$ is a color tree,
and the path $u_{i,6k-1}u_{i,6k+1}u_{i,6k-2}u_{i,6k}$ is a color tree. Suppose $x_i=1$. Then we give an anticlockwise coloring
such that $u_{i,6k}v_{i,2k}u_{i,6k+1}v_{i,2k-1}$ is a color tree ($u_{i,6m_i+1}=u_{i,1}$). If a variable $x_i$ is positive in clause $c_j$ for some $j\in[m]$,
then the path $w_{i,3k}u_{i,6k-2}u_{i,6k}w_{i,3k-2}w_{i,3k-1}u_{i,6k-1}$ is a color tree, and the path $u_{i,6k-3}u_{i,6k-5}u_{i,6k-2}u_{i,6k-4}$ is a color tree.
If a variable $x_i$ is negative in clause $c_j$ for some $j\in[m]$, then the path $u_{i,6k-1}w_{i,3k-1}u_{i,6k}u_{i,6k-2}$ is a color tree,
the path $w_{i,3k}u_{i,6k-3}w_{i,3k-1}w_{i,3k-2}$ is a color tree, and the path $u_{i,6k-3}u_{i,6k-5}u_{i,6k-2}u_{i,6k-4}$ is a color tree.
Without loss of generality, we can assume that the value of the first literal corresponding to the variable $x_i$ is true in clause $c_j$, see Figure \ref{mc4}.
The path $s_{j,2}ss_{j_1}s_{j,3}$ is a color tree. Notice that the edge $u_{i,6k-4}w_{i,3k}$ has been connected by some color tree.
Then we can construct a color tree on four vertices and a color tree on five vertices to connect the remaining edges in $I_j$ for each $j\in[m]$.
Now, all the nontrivial color trees are defined, and these color trees connect every edges in $\overline {G_\phi}$.
Then $f$ is a 1-color connection coloring of $G_\phi$. We have that $W_1(G)\leq7m-2+\sum_{i=1}^n 8m_i+\sum_{j=1}^m(3+2+2)=38m-2$.
Clearly, $cc_1(G_{\phi})>m(G_\phi)-38m+1$.

Suppose that $cc_1(G_{\phi})>m(G_\phi)-38m+1$. Let $f$ be an extremal 1-color connection coloring of $G_{\phi}$.
Then $f$ consists of several nontrivial color trees that connect all the edges in $\overline G_{\phi}$.
Clearly, the vertices of $V_0$ are in a same color tree; otherwise, we have to waste extra colors to connect the edges in $\{st: s,t\in V_0\}$.
Let $T$ be the color tree that contains the vertices of $V_0$. Suppose that the vertex $s_{j,4}$ is in $T$ for some $j\in[m]$.
Let $x_i$ be a variable corresponding to the literal $l_{j,1}$ in the clause $c_j$.
Then the vertex $w_{i,3k}$ is also in $T$ for some $k\in[m_i]$. Similar thing happens for the literals $l_{j,2}$ and $l_{j,3}$.
We have that $W_1(G)\geq 7l-2+\sum_{i=1}^n 9m_i+\sum_{j=1}^m(1+1+2)+1=38m-1$, a contradiction.
As a result, the vertex $s_{j,4}$ is not in $T$, and the color tree $T$ is exact $P$.
There are more than 5 wasted edges related to $s_{j,4}$ by Lemma \ref{p}.
We have that $W_1(G)\geq 7l-2+\sum_{i=1}^n 8m_i+\sum_{j=1}^m(5+2)=38m-2$.
Next we only need to consider the situation when equality holds since $cc_1(G_{\phi})>m(G_\phi)-38m+1$.
Then the variable garget $H_i$ has a ``clockwise" or ``anticlockwise" coloring consisting of some color trees $P_4$ or $P_6$ in Figure \ref{mc23};
otherwise, $f$ is not an extremal 1-color connection coloring. If $f$ is a clockwise coloring restricted to $H_i$ such that
the path $u_{i,6k-5}v_{i,2k}u_{i,6k-4}v_{i,2k-1}$ is a color tree, then $x_i=0$. If $f$ is an anticlockwise coloring restricted to $H_i$
such that the path $u_{i,6k-5}v_{i,2k}u_{i,6k-4}v_{i,2k-1}$ is a color tree, then $x_i=1$. Since the equality holds, there are only 5 wasted edges left
for the remaining edges in $I_j$ for each $j\in[m]$. Then there exists at least one literal in $c_j$ with true value for each $j\in[m]$.
Thus, $f$ is a 1-color connection coloring such that $\phi$ is satisfiable.

Finally, we have that $cc_1(G_{\phi})>m(G_\phi)-38m+1$ if and only if $\phi$ is satisfiable. The proof is now complete. \qed

\end{pf}

 Therefore, we have the following corollary.

\begin{cor}
Givin a graph $G$, it is NP-hard to compute the monochromatic connection number $mc(G)$ of $G$.
\end{cor}

\section{A good algorithm for the $k$-color connection number of trees}

Although to compute the $k$-color connection number for a general graph is NP-hard even for $k=1$,
it is polynomial-time to compute the number for some special graphs.
Let $T$ be a tree with $n$ vertices. Clearly, $cc_1(T)=1$. Next, we are going to find
an extremal $k$-color connecting coloring for trees with some nice properties, which will be used
for our algorithms.

\begin{lem}\label{treealg}
For a tree $T$, there exist an extremal $k$-color connection coloring $f$ of $T$ satisfying the following properties.
\begin{enumerate}
\item Each color forest $F_c$ is a color tree $T_c$.
\item If $T_c$ and $T_{c^\star}$ are two nontrivial color trees, then $T_c$ and $T_{c^\star}$ share a common vertex.
\item If $T_c$ is a nontrivial color tree and $v$ is a pendent vertex, then $cd\big(\{v\}, V(T_c)\big)\leq \lfloor\frac{k}{2}\rfloor$.
\end{enumerate}
\end{lem}

\begin{pf}
 Assume to the contrary, there exist a color forest $F_c$ containing two or more components.
 We do some contraction operations on $T$. First, we contract the edge $e\in E(T)$ if $E(F_c)$ lie in one component of $T\backslash e$.
 Next we contract all the edges of $F_c$. Finally, we pick a pendent vertex $v_0$ of the resulting graph,
 and let $e_0$ be the edge adjacent to $v_0$. Clearly, $v_0$ represents for a component $C_0$ of $F_c$.
 We recolor the edges of $C_0$ with color $f(e_0)$. Let $P=s\cdots t$ be the path between vertices $s$ and $t$.
 Then the new coloring $f^\star$ is still an extremal $k$-color connection coloring. However, the number of
 total number of components of $F_c$ decreases by one, a contradiction, and the first property is thus proved.

 For the second property, assume to the contrary, we have $cd\big(V(T_c),V(T_{c^\star})\big)\geq1$
 and the shortest path that uses minimum number of colors is $aa_0\cdots b$. Let $H$ be an induced subgraph by
 deleting the edges of the path $aa_0\cdots b$ from $E(T)$, and the component which contains $T_c$ is $A$,
 as well as the component which contains $T_{c^\star}$ is $B$. Without loss of generality,
 assume that for a pendant vertex $v\in V(A)$ we have $cd(v,a)\leq \lfloor \frac{k}{2}\rfloor$.
 Then we recolor an edge adjacent to $a$ in $T_c$ with color $f(aa_0)$. The new coloring $f^\star$ is still an
 extremal $k$-color connecting coloring. However, the number of edges of $T_c$ decreases by one, a contradiction.

 For the third property, again assume to the contrary, there exists a pendent vertex $v_0$ such that
 $cd\big(\{v_0\}, V(T_c)\big)\geq \lfloor\frac{k}{2}\rfloor+1$. Let $v_0\cdots t_0t$ be the path that uses $cd\big(\{v_0\}, V(T_c)\big)$ colors,
 and $v$ be a pendent vertex such that the path $v\cdots t_1t_2\cdots t$ goes through a pendent vertex $t_1$ of $T_c$.
 Then we exchange the colors of the edges $t_0t$ and $t_1t_2$. The new coloring $f^\star$ is still an extremal
 $k$-color connecting coloring. However, the length of the path that uses $cd\big(\{v_0\}, V(T_c)\big)$ colors decreases by one, a contradiction. \qed

\end{pf}

\begin{algorithm}
\caption{$cc_k(T)$ for $k$ is odd}
\label{odd}

\KwIn {A connected rooted tree $T$ with $V(T)=\{v_1, v_2, \cdots, v_n\}$; the root vertex $v_n$ lies on a path with length the diameter of
$T$ such that the distance between each end point of the path and $v_n$ is at least $\lfloor\frac{k}{2}\rfloor$; if $j > i$,
then $v_j$ is the parent of $v_i \in V(T)$; an odd integer $k$}
\KwOut { $a(v_n)+1$ }
\For{$i=1; i\leq n$}
    {
       $a(v_i)=0;
       b(v_i)=0$\;
    }
  \For{$i=1; i\leq n$}
    {
       \While{$v_j$ is the parent of $v_i$}
         {
           \eIf{$b(v_i)=\frac{k-1}{2}$}
           {
           $a(v_j)=a(v_i)+a(v_j)$\;
           $b(v_j)=b(v_i)$\;
           }
           {
           $a(v_j)=a(v_i)+a(v_j)+1$\;
           $b(v_j)=max\{b(v_j),b(v_i)+1\}$\;
           }
         }
    }
    return $a(v_n)+1$\;
\end{algorithm}

\begin{algorithm}[h]
\caption{$cc_k(T)$ for $k$ is even}
\label{even}

\KwIn {A connected rooted tree $T$ with $V(T)=\{v_1, v_2, \cdots, v_n\}$; the root vertex $v_n$ lies on the path with length the
diameter of $T$ such that the distance between each end point of the path and $v_n$ is at least $\lfloor\frac{k}{2}\rfloor$;
if $j > i$, then $v_j$ is the parent of $v_i \in V(T)$; an even integer $k$}
\KwOut { $a(v_n)+\Delta$}
\For{$i=1; i\leq n$}
    {
       $a(v_i)=0;
       b(v_i)=0$\;
    }
$S=\emptyset$\;
  \For{$i=1; i\leq n$}
    {
       \While{$v_j$ is the parent of $v_i$}
         {
           \eIf{$b(v_i)=\frac{k}{2}-1$}
           {
           $a(v_j)=a(v_i)+a(v_j)$\;
           $b(v_j)=b(v_i)$\;
           $S=S\cup \{v_iv_j\}$\;
           }
           {
           $a(v_j)=a(v_i)+a(v_j)+1$\;
           $b(v_j)=max\{b(v_j),b(v_i)+1\}$\;
           }
         }
    }
  $\Delta$ is the maximum degree of the subgraph induced by edge set $S$\;
  return $a(v_n)+\Delta$\;
\end{algorithm}

\begin{thm}\label{linear}
For a tree $T$ on $n$ vertices, Algorithms \ref{odd} and \ref{even} can compute $cc_k(T)$ in a linear time
of $n$.
\end{thm}

\begin{pf}
Let the root $u$ be one of the middle vertices (or vertex) on the path with length the diameter of the tree $T$.
We then use Breadth-First-Search algorithm to get a rooted tree $T$ with Property 1.
When $k$ is odd, the ``for" loops in Algorithm \ref{odd} computes the maximum number of distinct colors
outside the unique nontrivial color tree. Hence, Algorithm \ref{odd} computes $cc_k(T)$ by Lemma \ref{treealg}.
When $k$ is even, the ``for" loops in Algorithm \ref{odd} computes the maximum number of distinct colors outside
several combined nontrivial color trees. Hence, Algorithm \ref{even} computes $cc_k(T)$ by Lemma \ref{treealg}.
The running time of each iteration of the ``for" loops of Algorithm \ref{odd} or \ref{even} is $\mathcal{O}(1)$.
and so the running time is in linear of $n$. This completes the proof. \qed
\end{pf}

\end{document}